\documentclass[12pt]{article}
\usepackage{amstext,amsfonts,amssymb,amsthm,amsmath}
\usepackage{graphicx}
\usepackage[dvipsnames,fixpdftex]{xcolor}
\usepackage{pdfcolmk}
\usepackage[latin1]{inputenc}
\usepackage[margin=1.5cm,a4paper,twoside]{geometry}
\usepackage{hyperref}
\usepackage{setspace}
\usepackage{tikz}
\usepackage[percent]{overpic}
\usepackage[margin=10pt,font=small,labelfont=bf,format=plain,justification=justified]{caption}
\usepackage{wrapfig}

\hypersetup{colorlinks=false}
\hypersetup{pdfborder={0 0 0}}
\hypersetup{pdfpagemode=UseOutlines}
\hypersetup{pdfpagelayout=TwoColumnRight}

\newcommand{\vectdue}[2]{\left(\begin{matrix}#1\cr #2\end{matrix}\right)}

\renewcommand{\d}{\mathrm{\,d}}
\newcommand{\R}{\mathbb{R}}
\newcommand{\N}{\mathbb{N}}
\newcommand{\Z}{\mathbb{Z}}
\newcommand{\tr}{\mathrm{tr\,}}
\newcommand{\id}{\mathrm{id}}
\newcommand{\Id}{\mathrm{Id}}
\newcommand{\Div}{\mathrm{div\,}}

\newcommand{\particles}{\mathrm{part}}
\newcommand{\matr}{\mathrm{matr}}
\newcommand{\ext}{\mathrm{ext}}
\newcommand{\demag}{\mathrm{demag}}
\newcommand{\anis}{\mathrm{anis}}

\newcommand{\cell}{\square}

\newcommand{\unit}[1]{\ensuremath{\,\,\mathrm{#1}}}

\newcommand{\MPa}{\unit{M\,Pa}}
\newcommand{\GPa}{\unit{G\,Pa}}

\newcommand{\T}{\unit{T}}
\newcommand{\dist}{\mathrm{dist}}

\begin{document}
\begin{center}
  {\Large Hysteresis in Magnetic Shape Memory Composites:\\
         Modeling and Simulation}\\[2em]
         {\today}\\[2mm]
 {Sergio Conti$^1$, Martin Lenz$^2$, and Martin Rumpf$^{1,2}$}\\[2mm]
 {\em $^1$ Institut f\"ur Angewandte Mathematik,
Universit\"at Bonn, 53115 Bonn, Germany }\\
{\em $^2$ Institut f\"ur Numerische Simulation,
Universit\"at Bonn, 53115 Bonn, Germany }
\end{center}

\begin{abstract}
Magnetic shape memory alloys  are characterized by the coupling 
between a structural phase transition and magnetic one. This 
permits to control the shape change via an external magnetic
field, at least in single crystals. Composite materials with single-crystalline
particles embedded in a softer matrix have been proposed as a way to 
overcome the blocking of the transformation at grain boundaries.

We investigate hysteresis phenomena 
for small NiMnGa single crystals embedded in a polymer matrix
for slowly varying magnetic fields.
The evolution of the microstructure is studied within  the rate-independent variational framework proposed by Mielke and Theil (1999).
The underlying variational model incorporates linearized elasticity, micromagnetism, stray field
and a dissipation term proportional to the volume swept by the phase boundary.
The time discretization is based on an incremental minimization of the sum of energy and dissipation.
A backtracking approach is employed to approximately ensure the global minimality condition.

We illustrate and discuss the influence of the particle geometry (volume
fraction, shape, arrangement) and the polymer elastic parameters on the observed hysteresis
and compare with recent experimental results.
\end{abstract}

\section {Introduction}
Shape-memory alloys are crystalline materials which undergo a solid-solid phase transformation from a high-temperature,
high-symmetry austenitic phase to a low-temperature, low-symmetry martensitic one. The spontaneous shears induced
by the transformation are large but difficult to control; typical configurations after the phase transformation
consist of fine mixtures of 
different variants of the martensitic phase whose eigenstrains largely cancel each other, resulting in 
a very small average (macroscopic) net deformation.
Magnetic shape-memory (MSM) alloys are multiferroic materials, in the sense that besides the shape-memory effect they
are ferromagnetic. The different variants of the martensitic phase have different magnetic anisotropies, leading to a 
coupling between the magnetization and the eigenstrain. This permits to select one of the variants over the others,
and therefore to induce large eigenstrains, by the application
of external magnetic fields. The deformation strain in treated NiMnGa single crystals reaches $10\%$ \cite{Ullakko1996,Tickle99,Murray2000,Sozinov2002}.

Practical applicability of MSM single crystals for actuation and sensing is not easy, both because of the difficult production 
of single crystals and of their brittleness. In polycrystals, however, the transition is inhibited by grain boundaries. Indeed, if the
orientation of the grains is random then the 
eigenstrains of neighboring grains are typically not compatible with each other, this results in a blocking of the transition.
 Composite materials with single-crystalline
particles embedded in a softer matrix have been proposed as a way to 
overcome these difficulties \cite{Feuchtwanger2003,Hosoda2004,Feuchtwanger2005,Scheerbaum2007b,TianChen2009}, see
\cite{LiuScheerbaumKauffmannetal2012} for a recent review of the field.

Magnetic-field induced transformation has been demonstrated in composite materials, but 
the magnitude and even the presence of the effect depends strongly on many design parameters \cite{LiuScheerbaumKauffmannetal2012}. 
Typical examples are the stiffness of the polymer matrix, the size and shape of the particles and their density.
A theoretical investigation is therefore valuable not only to gain a better understanding of the material, but also 
to guide experimental search for the optimal material design. 
The static properties of MSM-polymer composites have been studied based on a variational model which couples micromagnetism and elasticity
both in the limit of small particles which contain no twin boundaries \cite{ContiLenzRumpf2007} 
and in the limit of large particles with a large number of twin boundaries \cite{ContiLenzRumpf2012}.  Both papers were
restricted to a static setting and  based on energy minimization without any account of time-dependent effects such as hysteresis. To the best of our knowledge both time-dependence and the physically relevant regime with particles of intermediate size 
have not been addressed yet.

We study here hysteresis in MSM-polymer composite materials, in the intermediate
regime in which few twin boundaries are present in each particle. We work in a quasistatic setting, using the 
rate-independent variational framework of Mielke and Theil \cite{MielkeTheil1999,MielkeTheilLevitas2002,MielkeTheil2004,Mielke2005} to account 
for the dissipation associated to the motion 
of interfaces. Our numerical results show how the hysteresis curves in composite materials depend on the 
material parameters and geometry and, we expect, will be valuable as guides for subsequent material
development.

\begin{figure}
\begin{center}
\begin{overpic}[width=0.49\linewidth]{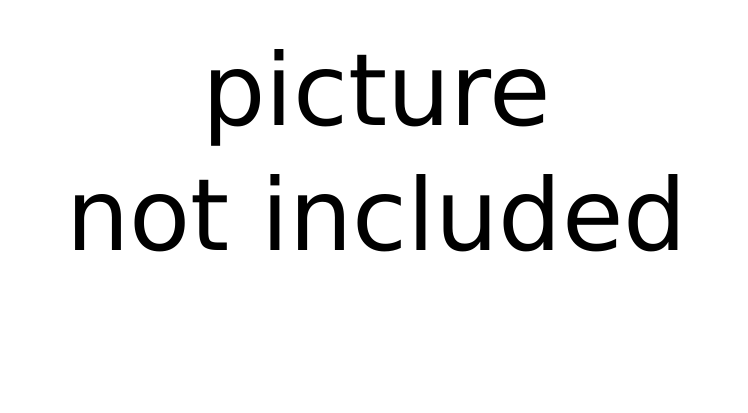}
 \put (3.5,3.5) {\color{white}\rule{4mm}{4mm}}
\end{overpic}%
\hfill
\includegraphics[width=0.49\linewidth]{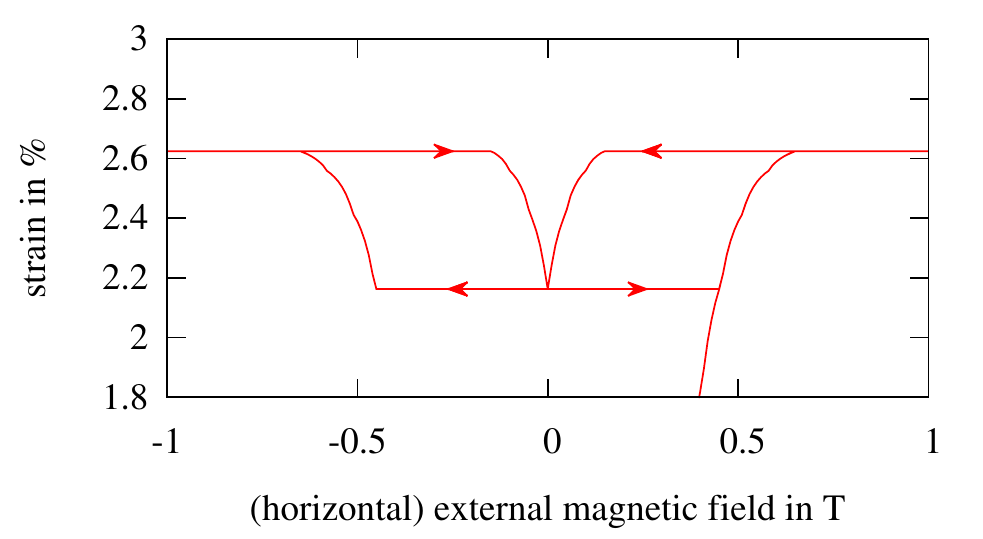}
\end{center}
\caption{Hysteresis in NiMnGa single crystals (left, 
reprinted from  \cite[Fig.~3a]{StrakaHeczko2005}, with permission
from Elsevier) 
compared to simulation of single particle in composite (right). In both cases the material is subject to an external magnetic field, which is 
increased from 0 up to 1 Tesla and back to 0, then and up to 1 Tesla in the opposite direction and back again. In case of the NiMnGa single crystal, a constant compressive stress ($1 \MPa$) pushes 
the deformation back towards the reference configuration. 
The simulation is done with $E=1 \MPa$ for the polymer.}
\label{fig:experiment_simulation}
\end{figure}

The starting point of our considerations is a static model for magnetic shape memory materials 
which couples micromagnetism and elasticity \cite{ContiLenzRumpf2007,ContiLenzRumpf2008}. The phase transformation enters the model via a phase index function $p$, 
which characterizes the local phase.  The eigenstrain entering the elastic energy density and the anisotropy of the
magnetic energy then depends on this phase index, since different martensite variants have different magnetic easy axes 
and different transformation strains, details are discussed in Section \ref{secmodel} below.
With changing external parameters  (magnetic field or mechanical pressure) the interface between different variants
changes. The motion of the twin boundary is, however, coupled to a dissipation. This makes the problem history-dependent and
generates hysteresis loops.
Hysteresis can be modeled without describing the fast timescales of elastic and magnetic oscillations, if
the external forces are 
changed only slowly. Mielke and Theil \cite{MielkeTheil1999,Mielke2005} 
proposed a rate independent modeling framework where one assumes that the energy from the fast oscillations is completely dissipated. The amount of energy dissipated depends on the path of the 
transformation. 
The model can then be formulated 
solely based on the dissipation and the (elasto-magnetic) energy.
Rate-independent models along these lines have been studied both for shape-memory alloys \cite{AuricchioMielkeStefanelli2008}
and for magnetic shape-memory alloys \cite{BessoudStefanelli2011}, without resolving the spatial details of the microstructure.

Figure~\ref{fig:experiment_simulation} shows a measured hysteresis loop in a 
magnetic shape memory single crystal compared to a simulation based on the rate-independent model for a composite.
The qualitative aspect of the hysteresis loop is very similar, although the magnitude of the strain change is significantly
reduced in the composite. The characteristics of this hysteresis loop are discussed in Section \ref{secresults} below in detail.
The simulation results we present here are, to the best of our knowledge, the first detailed description
of hysteresis in composites, since quantitative hysteresis measurements are quite difficult 
\cite{KauffmannScheerbaumLiuetal2012,LiuScheerbaumKauffmannetal2012}.

\medskip

\section {Model}\label{secmodel}

\begin{figure}
\begin{center}
\colorlet{colvariantA}{cyan!40!white}
\colorlet{colvariantB}{blue!50!white}
\colorlet{colmartensite}{colvariantA!50!colvariantB}

\begin{tikzpicture}[scale=0.6]
\begin{scope}[shift={(-6.2,0)}]
  \draw (-2.5,-4) node {(a)};
  \fill[yellow] (-3,-3) rectangle (3,3);
  \foreach \y in {-2.75, -2.25, ..., 2.8} {
  \foreach \x in {-2.75, -2.25, ..., 2.8} {
  \fill[colmartensite] (\x,\y) circle (.1666);} }
  \draw (0,0) rectangle (0.5,0.5);
\end{scope}
\begin{scope}
  \draw (-2.5,-4) node {(b)};
  \fill[yellow] (-3,-3) rectangle (3,3);
  \fill[colmartensite] (0,0) circle (2);
  \draw (-2,-2) node {$\Omega$};
  \draw (0,0) node {$\omega$};
\end{scope}
\begin{scope}[shift={(6.2,-0.3)}]
  \draw (-2.2,-3.7) node {(c)};
  \fill[yellow] (-2.7,-2.7) rectangle (3.3,3.3);
  \begin{scope}
  \clip (-3,-3) -- (3,3) -- (-3,3);
  \fill[colvariantA] (0,0) ellipse (1.7 and 2.3);
  \end{scope}
  \begin{scope}
  \clip (-3,-3) -- (3,3) -- (3,-3);
  \fill[colvariantB] (0,0) ellipse (2.3 and 1.7);
  \end{scope}
  \draw (-0.4,1.0) node {$p=1$};
  \draw (0.5,-0.9) node {$p=2$};
\end{scope}
\begin{scope}[shift={(12.4,-0.1)}]
  \draw (-2.4,-3.9) node {(d)};
  \fill[yellow] (-2.4,-3.4) rectangle (2.6,3.6);
  \begin{scope}
  \clip (-3,-4) -- (4,3) -- (-3,3);
  \fill[colvariantA] (0,0) ellipse (1.7 and 2.3);
  \end{scope}
  \begin{scope}
  \clip (-3,-4) -- (4,3) -- (3,-3);
  \fill[colvariantB] (-0.3,-0.3) ellipse (2.3 and 1.7);
  \end{scope}
  \draw (-0.2,0.8) node {$p=1$};
  \draw (.85,-1) node {$p=2$};
  \draw[thick,->] (-1.8,-3.9) -- (2,-3.9) node[right] {$H$};
\end{scope}
\end{tikzpicture}
\end{center}
\caption{Sketch of the configuration in one periodic cell
containing a single magnetic shape memory particle in a polymer matrix. 
(a) The sample contains many small particles, which we assume for simplicity to be all equal
and periodically distributed. In (b-d) only one unit cell is shown.
(b) The reference configuration is the ground state of the austenite phase,
with no eigenstrain. (c) In the martensite phase different variants
with different eigenstrains coexist in the particle, separated
by a small number of twin boundaries (here only one). By elastic compatibility
the twin boundaries are straight. 
(d) When applying an external magnetic field, variant $1$ is preferred due to its horizontal 
anisotropy, so a part of variant $2$ transforms to variant $1$ (i.e. the twin boundary moves). Variant $1$ is shorter in the horizontal direction, which leads to a deformation of the whole composite.}
\label{fig:configuration}
\end{figure}
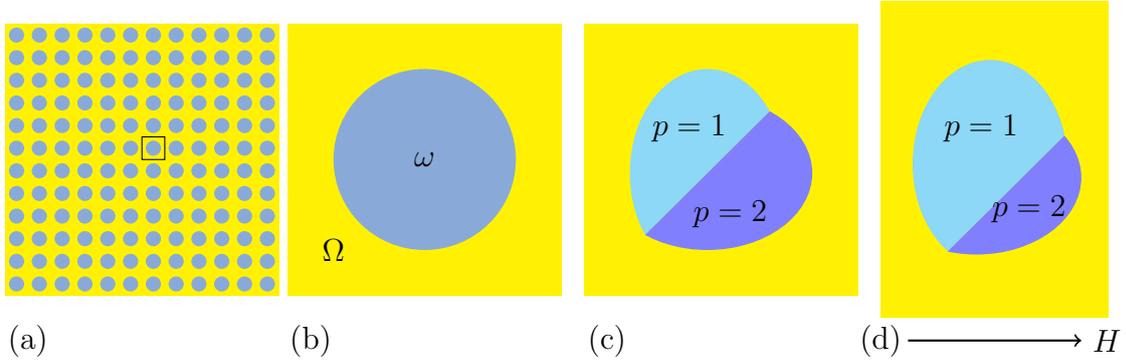

\paragraph {Energy of an MSM--polymer composite.} Our numerical implementation is restricted to two dimensions, already providing both qualitative and quantitative insight in the underlying hystereses.
Hence, we restrict here to a description of the model in 2D. The extension of the model to 3D is straightforward. 
We consider MSM particles $\omega_i$, $i=1,2,\ldots N$, and a polymer matrix $\Omega\setminus\omega$ with $\omega=\bigcup_{i=1}^N\omega_i$, see Figure
\ref{fig:configuration}(a) for an illustration. The  external magnetic field $H:[0,T]\to \R^2$  
is spatially uniform but time dependent. We model the deformation of the particles and the polymer with (linearized) elasticity and denote by 
$u:[0,T]\times\Omega\to \R^2$
the elastic displacement and by $\epsilon[u] = \frac12 \big(\nabla u + (\nabla u)^T\big)$ the linearized strain. 
The magnetic behavior of the MSM material is described by a magnetization field $m:[0,T]\times\R^2\to\R^2$. The phase index $p:[0,T]\times\omega\to\{1,2\}$
couples transformation strain and magnetic anisotropy. 
The jump set of $p(t,\cdot)$ (i.e. the twin boundaries) is denoted by $J_p(t,\cdot)$.
For the elasto-magnetic energy, we consider the configuration $u(t),m(t),p(t)$ 
at a specific time $t$ and for a fixed external magnetic field $H(t)$ 
(with a slight misuse of notation, we write $u,m$ and $p$ also for the formal parameters
of the energy).

We assume the polymer to be an isotropically elastic material, which is strain free in the reference configuration. The elasticity tensor of the particles is assumed to have a cubic symmetry, the 
phase index selects the martensitic variants and their eigenstrains. 
The reference configuration is usually assumed to be the austenitic phase, with respect to which the eigenstrains of the two variants are
\begin{equation}\label{eqdefeps12}
  \epsilon_1=\vectdue{-\epsilon_0 & 0}{0 & \epsilon_0} \hskip3mm\mbox { and }\hskip3mm
  \epsilon_2=\vectdue{\epsilon_0 & 0}{0 & -\epsilon_0}.
\end{equation}
In most of our simulations we assume the polymer to be solidified with the MSM particles in the austenitic phase, so that the
eigenstrain of the polymer is zero. Of course it is also
possible to model the situation in which the polymer is solidified with the particles already in one of the two martensitic
variants. Then it is more appropriate to take that state as reference, and $\epsilon_1-\epsilon_2$ and $0$ as eigenstrains of the two phases.

The orientation of the strains in (\ref{eqdefeps12}) is relative to the crystal lattice of the MSM alloy. The actual lattice orientation of the particles in the composite 
is described by a rotation $Q:\omega\to SO(2)$, which we suppose to be 
constant in each particle.
We use a subscripts $\Omega$ and $\omega$, respectively, to indicate the relevant domain of integration for the energy.
The elastic energy of the particles is given by
 \begin{equation}\label{eqelastparticles}
   E_{\omega}^\particles[u,p]=\int_\omega
  W^\particles((\epsilon [u](x))Q(x)-\epsilon_{p(x)})\d x
 \end{equation}
with an elastic energy density 
\begin{equation*}
W^\particles(\epsilon)= \tfrac12 C_\particles \epsilon : \epsilon = 
\frac12 C_{11} (\tr \epsilon)^2 + (C_{12}-C_{11}) \epsilon_{11}
\epsilon_{22} + 2 C_{44} \epsilon_{12}^2\
\end{equation*}
parametrized by the cubic elastic constants $C_{11}$, $C_{12}$ and $C_{44}$.
The elastic energy of the polymer matrix is defined as 
\begin{equation}\label{eqelastmatrix}
  E_{\Omega \setminus \omega}^\matr[u]=\int_{\Omega\setminus\omega} W^\matr(\epsilon [u](x))\d x + \int_{\Gamma_N} u \cdot g \d \mathcal{H}^1\,,
\end{equation}
where $W^\matr(\epsilon)= \tfrac12 C_\matr \epsilon : \epsilon = \tfrac12 \lambda (\tr\epsilon)^2 + \mu |\epsilon|^2$  with isotropic elastic constants $\lambda$ and $\mu$. 
Here $\Gamma_N\subset\partial\Omega$ is the Neumann boundary, and $g$ represents the applied surface traction. Dirichlet data on another subset
$\Gamma_D\subset\partial\Omega$ can then be imposed by restricting the set of admissible displacements $u$.

Micromagnetism has to be considered in the physical, i.e., the deformed configuration.
Coordinates in the deformed configuration are denoted by $y$ 
in contrast to $x$ for the reference configuration. 
The deformation itself is denoted by $v(x)=x+u(x)=(\id+u)(x)$.
We scale the magnetization by the saturation magnetization $M_s$ of the MSM alloy.
Thus $|m|=1$ on $v(\omega)$, and $m=0$ otherwise.

The micromagnetic energy 
describes the interaction with the (time dependent) external magnetic field (Zeeman energy $E_\omega^\ext$), 
the demagnetization (stray) field (magnetic exchange energy $E_{\R^2}^\demag$)
and the magnetic anisotropy ($E_{\omega}^\anis$). The anisotropy depends on the martensitic variant encoded by the phase index $p$ and we assume the magnetic easy axis to coincide with the 
shortening axis of the eigenstrain. We define
\begin{align}
  E_\omega^\ext[H,u,m]& =-\frac{M_s}{\mu_0} \int_{v(\omega)} H \cdot m\d y\,,\\
  E_{\R^2}^\demag[m]&= \frac{M_s^2}{\mu_0} \int_{\R^2} \frac12  |\nabla \psi|^2 \d y\,,\\
  & \qquad \mbox{ where } \Delta\psi = \Div m \mbox { distributionally in $\R^2$, } \\
  E_{\omega}^\anis[u,m,p]&=K_u\int_{v(\omega)} \phi_{p(v^{-1}(y))} \left( (R_{\nabla v \circ v^{-1}} Q(v^{-1}(y)))^T m \right) \d y\,,\\
         & \qquad \mbox{ where } \phi_2(m)=m_1^2 \,, \phi_1(m)=m_2^2.
\end{align}
Here, $\mu_0$ is the permeability of vacuum, $K_u$ an anisotropy constant (an energy density, since $m$ is dimensionless), and  $R_{\nabla v \circ v^{-1}}\in SO(2)$ the rotational part of $\nabla v$ 
at $v^{-1}(y)$. In all simulations we use the approximation $R \approx \Id+ \frac {1}{2} (\nabla v-(\nabla v)^T)$. 
The stray field potential $\psi$ is defined on $\R^2$ and has zero boundary data at infinity.

In our simulations we focus on the situation that inside each particle only one twin boundary is present. Due to elastic
compatibility they are then approximately straight. In particular, the phase index $p$ is piecewise constant inside each particle, 
as illustrated in Figure \ref{fig:configuration}(c). The exchange term can then be replaced by a term 
penalizing the length of the interfaces, therefore we do not include it explicitly in the model.

The total energy depending on the state $(u,m,p)$ and explicitly on time $t$ (based on the time varying external field) is given by
\begin{equation} 
E[t,u,m,p] = E_{\Omega \setminus \omega}^\matr [u] + E_{\omega}^\particles [u,p] + E_\omega^\ext[H(t),u,m] + E_{\R^2}^\demag[m] 
+ E_{\omega}^\anis [u,m,p]\,. 
\end{equation}
\bigskip

\paragraph {Rate-independent evolution model.} 
We study evolution in the rate-independent limit, i.e. the applied fields evolve on a timescale which is much slower than the
one of the internal processes and of internal equilibration. Hence, velocities play no role and the configuration solely depends on the path 
in state space and not on the rate at which the state changes along the path. Correspondingly the evolution of the system depends
only on the path taken by the external forces and fields, and not by their rate of change.

The key assumption in the models of \cite{MielkeTheil1999,Mielke2005} is that at any
moment in time the system is in the state which is energetically most favorable, where one has to correctly account for the
dissipation which would be incurred to move to a different state. Precisely, 
if a system is initially in a state $S_0$, then the state
$S$ at time $t$ is characterized as a minimizer of $E(S,t)+\mathrm{Diss}(S_0,S)$, where $E(S,t)$ is the energy of 
state $S$ at time $t$ (with the appropriate
external forces and boundary conditions) and $\mathrm{Diss}(S_0,S)$ the dissipation associated to the transition from $S_0$ to $S$ along the path in state space.

In our setting the state is described by the triple $(u,m,p)$. Whereas elasticity is typically considered to be a reversible phenomenon,
which generates little dissipation, the phase transition has a large dissipation. This can be made quantitative for example by the size
of the hysteresis loop in a single crystal illustrated in Figure \ref{fig:experiment_simulation}.  Therefore we assume the dissipation 
to depend only on the phase index $p$.
 
We assume that the dissipated energy is proportional to the transformed volume with a proportionality constant $\kappa$, and define a dissipation distance $D[p,q]$ between two different 
configurations ($p$ and $q$ being the respective phase indices)
as
\begin{equation}
 D[p,q]  = \int_\Omega \kappa |q(x)-p(x)| \d x\,.
\end{equation}
The corresponding accumulated dissipation $\mathrm{Diss}[\cdot,[t_1,t_2]]$, which depends on the entire evolution of the
phase index $p$ for times from $t_1$ to $t_2$, is defined by subdividing the interval $[t_1,t_2]$ into subintervals, and then taking
the maximum over all possible decompositions,
\begin{equation}
 \mathrm{Diss}[p;[0,t]]  = \sup \{ \sum_{i=1}^n D[p(t_{i-1}),p(t_i)] \,|\, n \in \N, 0 \leq t_0 < t_1 < \cdots < t_n \leq t \} 
\end{equation}
(a similar procedure gives the usual definition of dissipation length of a curve in state space).

A solution to the rate independent evolution problem \cite{MielkeTheil1999,Mielke2005} can now be defined solely 
in terms of the energy and the dissipation. Given the external field $h(t)$, the solution should verify at all times $t$ the 
stability condition
\begin{equation}
\tag{S}\label{eq:S}  E[t,u(t),m(t),p(t)]\leq E[t,\widetilde{u},\widetilde{m},\widetilde{p}]+D[p(t),\widetilde{p}]\quad\mbox{for all }\widetilde{u},\widetilde{m},\widetilde{p}\,, 
\end{equation}
and the energy balance
\begin{equation}
\tag{E}\label{eq:E}  E[t,u(t),m(t),p(t)]+\mathrm{Diss}[p;[0,t]]-E[0,u(0),m(0),p(0)]=\int^{t}_{0}{\tfrac{\partial}{\partial t}} E[s,u(s),m(s),p(s)] \d s\,,
\end{equation}
where $\tfrac{\partial}{\partial t}E$ represents the partial derivative of the energy with respect to time.
\paragraph{Physical constants.}
We choose parameters that match the experimentally known values for NiMnGa particles,
as in \cite{ContiLenzRumpf2007}.
In  the magnetic energy we use $\frac{M_s}{\mu_0} = 0.50 \frac{\MPa}{\T}$, $\frac{M_s^2}{\mu_0} = 0.31 \MPa$
\cite{Murray2000,LikhachevUllakko2000} and 
$K_u= 0.13 \MPa$ \cite{OHandley1998,Soederberg2005}.
The elastic constants used for NiMnGa are $\epsilon_0 = 0.058$, $C_{11} = 160 \GPa$, $C_{44}= 40 \GPa$, $C_{11}-C_{12} = 4 \GPa$ \cite{StipcichPlanes2004,ZhaoDaiCullenWuttig2007,HeczkoKopecekMajtasLanda2011}.
For the polymer the elastic modulus is $E= 1\MPa$ unless specified otherwise (the polyurethane in \cite{KauffmannScheerbaumLiuetal2012} has $E=2\MPa$) and the Poisson ratio $\nu= 0.45$. 

\paragraph{Periodic composite structures.} 
In most practically relevant cases, the number of particles in the composite workpiece is very large. Instead of 
simulating the complete composite  with its very complicated geometric structure 
we strive for a homogenization approach which asks for the 
effective macroscopic material properties based on the solution of suitable microscopic problems \cite{Willis1981,Oleinik1992,Braides98,DonatoCioranescu,Milton2002}.
To this end,  we study periodic configurations, where the particle geometry in a single cell is simple, as illustrated in Figure \ref{fig:configuration}(a).
In a static context based on energy minimization this procedure was already discussed heuristically 
in \cite{ContiLenzRumpf2007}, a homogenization result was then proven rigorously in  \cite{ThesisPawelczyk14}.

We assume that the microstructure is obtained by downscaling a 
fixed structure defined on the unit square $\Omega^\cell = [0,1]^2$ and assume for simplicity that 
the MSM phase $\omega^\cell$ on the reference cell does not intersect the boundary. For a composite workpiece $\Omega$ with this fine scale structure we obtain
the periodic MSM phase $\omega = \Omega \cap \varepsilon (\omega^\cell + \Z^2)$, with
$\varepsilon>0$ a small parameter describing the periodicity of the material.
The associated lattice orientation is given by 
$Q(x) = Q^\cell(\frac{x}{\varepsilon} \mod \Z^2)$, where $Q^\cell$ is the lattice orientation on the reference cell.
Now, the theory of homogenization separates the microscopic scale, which
resolves the full complexity of the given microstructure, from the macroscopic scale, which 
describes the effective behavior of the composite workpiece in the limit $\varepsilon \to 0$.
Thus, we restrict ourselves here to the reference configuration $\Omega^\cell$ with particle domain $\omega^\cell$ 
and study the energy and dissipation functional on this reference cell to compute the effective hysteresis properties. 
Furthermore, we restrict here to macroscopic affine displacements $x\to A x$ with a symmetric matrix 
$A \in \R^{2\times 2}$ and assume correspondingly
\begin{equation}\label{eqbcaffineper}
u(x+e_i) = u(x) + A e_i \,,\hskip1cm i=1,2
\end{equation}
 on the boundary of the reference domain $\Omega^\cell$. Let us emphasize that we 
investigate only linearized elasticity and assume invariance with respect to infinitesimal 
rotations. 

If we consider a composite material covering $\R^2$ ($\Omega = \R^2$) the associated stray field $\psi_{\R^2}$ is periodic.
For general workpieces $\Omega \subsetneq \R^2$ the resulting stray field $\psi_\Omega$ can be split into the periodic component $\psi_{\R^2}$ 
and a non periodic component $\psi_\Omega - \psi_{\R^2}$. For a detailed discussion of this we refer to \cite{ContiLenzRumpf2007}.
In most of the computations in this paper we treat only the periodic component $\psi_{\R^2}$. In the example in Fig. \ref{fig:macrostrayfield} below
we investigate the full model for a circular workpiece, where the effective stray field correction $\psi_\Omega - \psi_{\R^2}$ can be explicitly computed
via a modification of the external field $H$ (cf. \cite{ContiLenzRumpf2007}).

The energy is not convex. Thus, minimizers of the energy do not necessarily share the same periodicity as the energy itself, 
in particular minimizing the energy over periodic configurations with period larger than one 
might lead to lower energy values (cf. \cite[Section 6.3]{ContiLenzRumpf2007}). In this paper 
we do not address this effect and always work with the smallest possible period. For details of the homogenization procedure
we refer to \cite{ThesisPawelczyk14}.

\section {Time Discretization}

The formulation of the rate-independent evolution problem given in \eqref{eq:E} and \eqref{eq:S} lends itself
naturally to time discretization \cite{MielkeTheil1999,Mielke2005}.
We fix time steps $t_0, t_1, t_2 , \ldots$ and ask for approximations 
$u_i \approx u(t_i), m_i \approx m(t_i), p_i \approx p(t_i)$.
The natural time discretization is then given by the time-incremental formulation
\begin{displaymath} 
\tag{T}\label{eq:T} (u_i, m_i, p_i) \quad \mbox { minimizes } \quad
E[t_i,\widetilde u, \widetilde m, \widetilde p] + D [p_{i-1},\widetilde p]\,.
\end{displaymath}
This means that at each of the time steps the state of the material minimizes the sum of the 
energy (which depends explicitly on time through the forcing terms and surface tractions) and the dissipation,
taken with respect to the state at the previous time.
Taking into account the triangle inequality for the dissipation distance $D[\cdot, \cdot]$ it immediately follows that a solutions of \eqref{eq:T} also satisfy \eqref{eq:S}. 
The energy balance \eqref{eq:E} is fulfilled approximately. Indeed, 
using the material state at time $t_{i-1}$ as a comparison function in \eqref{eq:T} and taking into account $D[p_{i-1},p_{i-1}]=0$,
implies $E[t_i,u_i,m_i,p_i]+D[p_{i-1},p_i] \le E[t_{i},u_{i-1},m_{i-1},p_{i-1}]$. Now one obtains
 \begin{equation}
 \tag{E--}\label{eq:EL} E[t_i,u_i,m_i,p_i]+D[p_{i-1},p_i]-E[t_{i-1},u_{i-1},m_{i-1},p_{i-1}]\leq \int^{t_i}_{t_{i-1}}{\tfrac{\partial}{\partial t}} E[s,u_{i-1},m_{i-1},p_{i-1}] \, ds\,.
 \end{equation}
Since the material state at time $t_{i-1}$ is stable, using
\eqref{eq:S} at time $t_{i-1}$ with $(u_i,m_i,p_i)$ as a comparison state leads to
$E[t_{i-1},u_{i-1},m_{i-1},p_{i-1}] \le
E[t_{i-1},u_i,m_i,p_i]+D[p_{i-1},p_{i}]$. As above, we rewrite the inequality in the form
 \begin{equation}
\tag{E+}\label{eq:ER} E[t_i,u_i,m_i,p_i]+D[p_{i-1},p_i]-E[t_{i-1},u_{i-1},m_{i-1},p_{i-1}]\geq \int^{t_i}_{t_{i-1}}{\tfrac{\partial}{\partial t}}E[s,u_{i},m_{i},p_{i}] \, ds\,.
 \end{equation}
The two conditions \eqref{eq:EL} and \eqref{eq:ER} constitute an approximate, local in time version of the energy condition
\eqref{eq:E}.
Our numerical solution is based on solving \eqref{eq:T} on a suitably chosen time discretization. Since 
\eqref{eq:T} demands for a global minimum, but numerical algorithms normally only locate local minima, 
a backtracking scheme is employed,  see discussion in Section \ref{eqc:timediscr} below.
 
\section {Space Discretization} \label{spacediscrete}
The space discretization is based on a direct boundary element ansatz using a collocation discretization 
with piecewise constant (for the demagnetization field) and piecewise affine (for the displacement) ansatz functions \cite{At97,ClRi78}. 
To this end, the boundary of the computational domain (for a single cell the unit square) and the 
particle--matrix--interface are approximated by polygons. 
The twin boundaries also have to be discretized. 
This results in a partitioning of the domain in the polymer matrix and two subdomains (the twins) per particle (cf. Fig. \ref{fig:configuration}). 
Here we restrict to straight lines as twin boundaries within the MSM particles. In the time discrete model a system of partial differential equations has to be solved for
the elasticity and for the demagnetization field. The equations are linear and the corresponding coefficients are constant on each component of the partition. 
Thus a boundary element discretization can be set up on each part separately, coupled via appropriate interface conditions on the different types of boundary: boundary of the periodic cell, 
particle-matrix-interface, twin boundary. Special considerations are necessary at the triple point where the twin boundary meets the particle boundary, details are described below in the discussion 
of the treatment of the different energy contributions.

\paragraph{Degrees of freedom.} Due to our restriction to straight twin boundaries inside 
each particle, these interfaces can be described by two parameters (angle and distance from the particle 
center). In the application, the MSM particles are much harder than the polymer matrix. Thus, we may assume that the twins within each particle undergo just affine deformations. In case of two twins 
per particles, we have 12 degrees of freedom for the two affine deformations.
The rank-1-condition along the interface leaves only 8 degrees of freedom: Shift and rotation of the whole particle and tangential stretching  aligned with the twin boundary ($2+1+1$), separate tangential shear 
and normal stretch on both sides ($2\cdot(1+1)$). The eigenstrain can actually be realized by a tangential shear of opposite sign on both sides of the interface.
The affine transformation $x\mapsto x + Ax$ of the reference cell is described by $3$ degrees of freedom (due to the symmetry of $A$).  In addition, 
we assume the magnetic domain walls to coincide with the twin boundaries. In particular, on each of these domains the magnetization is described by one degree of freedom, namely the direction of the 
magnetization.
In our model, only the polymer elasticity and the demagnetization field are actually discretized using a collocation boundary element method with degrees of freedom associated to the vertices of the polygonal MSM-polymer interface $\partial \omega^\cell$
and vertices on the boundary $\partial \Omega^\cell$ of the reference cell.

\paragraph{Particle elasticity, anisotropy, interaction with magnetic field.}
Since the deformation of each twin is assumed to be an affine function, the elastic energy of a particle depends on the set of associated piecewise constant strains and can  be easily 
evaluated. The anisotropy and the Zeeman energy are computed on the deformed domain. Thus, they depend on the magnetization and the affine deformations on the twins and the evaluation is immediate. 

\paragraph{Polymer elasticity.}
The elasticity of the polymer has to be resolved in more detail. It has to accommodate for the different deformations of the twinned particles and the motion of the twin boundary, and so the strain 
varies significantly throughout the matrix.  We assume that the polymer relaxes instantaneously depending on the deformation of the particles and the macroscopic strain. Thus, the associated (static) 
linearized elasticity problem has to be solved on $\Omega^\cell \setminus \omega^\cell$ with two different types of boundary conditions, namely Dirichlet boundary conditions on $\partial\omega^\cell$ 
depending on the (given) deformation $u_p$ of the particle and periodic boundary conditions with the affine offset $x\mapsto A x$ on the boundary of the unit cell $\Omega^\cell = [0,1]^2$. 
We obtain
\begin{align*}
- \Div \Big( \lambda \tr\epsilon[u](x)) \id + 2 \mu \epsilon[u](x) \Big) = 0 & \mbox { for } x \in \Omega^\cell \setminus\omega, \\
u(x) = u_p(x) & \mbox { for } x \in \partial\omega,\\
u(x+e_i) = u(x) + A \, e_i & \mbox { for } x \in \partial \Omega^\cell, \ i=1,2. 
\end{align*}
Once displacement and normal stresses are computed on $\partial \Omega^\cell \cup \partial \omega^\cell$ by the boundary element method 
the elastic energy can be evaluated by a straightforward integration by parts \cite{ContiLenzRumpf2007}.

When using the elastic energy in the descent algorithm  additional regularizations are required: The displacement on the particle boundary is given by the affine deformation of the two twins, and enters the boundary element computation of the polymer elasticity as a boundary condition. Usually this boundary displacement will have a kink where the twin boundary meets the particle boundary. 
The kink is actually caused by the jump of the tangential shear (dominant) and the normal stretch at the twin boundary. The associated component of the deformations, denoted by  $u^{\mathrm{kink}}_p$ depend linearly on the distance $\dist(x)$ to the plane of the twin boundary, i.e. we write $u^{\mathrm{kink}}_p=u^{\mathrm{kink}}_p(\dist(x))$.
In our implementation, we take into account the smooth approximation $\dist^\delta(x) =  (\dist^2(x)+\delta^2)^{\frac12}-\delta$ for the distance $\dist(x)$ 
and replace $u^{\mathrm{kink}}_p(\dist(x))$ by $u^{\mathrm{kink}}_p(\dist^\delta(x))$.
The parameter $\delta$ has to be chosen of the order of the grid size $h$ used to discretize $\partial \omega^\cell$. In fact, we choose $\delta=2h$.

\paragraph{Demagnetization field.} 
The stray field $\psi$ 
solves  $\Delta\psi = \Div m$ distributionally on $(\id + v)(\Omega^\cell)$ with periodic boundary conditions. Since $m$ is piecewise constant, this can be expanded to 
\begin{align} 
\Delta \psi &= 0 & \mbox { in } v(\Omega^\cell) \setminus \partial \omega \cup J_p), \nonumber \\
[\nabla \psi \cdot \nu] &= [m \cdot \nu] & \mbox { on } v(\partial \omega \cup J_p), \label{eq:strayjump} \\
\psi(v(x+ e_i)) &= \psi(v(x)) & \mbox { for } x \in \Omega^\cell, \ i=1,2. \nonumber
\end{align}
Here, $[f]$ indicates the jump of the function $f$ at the interface along the direction of the normal $\nu$.
Now, we consider a splitting $\psi = \psi_J + \psi_P$, where $\psi_J$ is defined on whole $\R^2$ with $\Delta \psi_J = 0$ a part from the jump set $J_p$ and fulfills the jump condition \eqref{eq:strayjump}. Furthermore, we assume that $\psi_P$  solves $\Delta \psi_P = 0$ on $v(\Omega^\cell)$ without  jump.
The boundary conditions for $\psi_P$ follow from the periodic boundary conditions for $\psi$ and the splitting assumption, i.e.
\[\psi_P(v(x+ e_i)) = \psi_P(v(x)) + \psi_J(v(x)) - \psi_J(v(x+ e_i)) \mbox { for } x \in \partial \Omega^\cell, \ i=1,2. \]
The component $\psi_J$ can then be directly computed based on the integral representation 
\[ \psi_J(x) = \int_{v(\partial \omega \cup J_p)} G(x-y) \, [m\cdot \nu](y) \d y, \]
where $G(x)= -\frac1{2\pi}\log|x|$ is the fundamental solution of the Laplacian.
Using integration by parts one obtains for the component $\psi_J$ of the demagnetization energy 
\begin{align*}
\int_{v(\Omega^\cell)} \frac12 |\nabla \psi_J|^2 \d y & = \int_{v(\partial \omega^\cell \cup J_p)} \psi_J(x) \, [m\cdot \nu](x)  \d x \\
& =  \int_{v(\partial \omega^\cell \cup J_p)} \int_{v(\partial \omega^\cell \cup J_p)}  G(x-y) \, [m\cdot \nu](y)  \, [m\cdot \nu](x)  \d x \d y\,.
\end{align*}
For $\psi_P$ we employ the same boundary element strategy as for the polymer elasticity discussed above. The remaining integrals of $|\nabla \psi_P|^2$ and $\nabla\psi_P \cdot \nabla\psi_J$ for the complete  demagnetization energy can analogously be rewritten als integrals over the interfaces involving values of $\psi_P$ and $\psi_J$ on these interfaces and the normal jump of $m$.
The above regularization at the triple points where twin boundary and particle boundary meet 
can be omitted since the singular component of the stray field that is affected by the triple point is computed exactly.

\section {Implementation of the Time Discretization} 
\label{eqc:timediscr}
\paragraph{Energy descent.} 
States $(u,m,p)$ are described in the spatially discrete case by a vector $z \in \R^N$, 
where $N$ is the total number of degrees of freedom; i.e. we write $(u(z),m(z),p(z))$. Let us suppose that an index set $I_p \subset \{1,\ldots, N\}$ identifies the degrees of freedom describing the phase $p$. We use a gradient descent scheme for the numerical solution of the minimization problem \eqref{eq:T} in each time step. To this end a descent direction has to be computed. Due to the above described regularization the energy $E$ is differentiable but the dissipation distance $D$ is not. Indeed, $p \mapsto D[p_{i-1},p]$ possesses only a subgradient at $p=p_{i-1}$. 
Thus, we proceed as follows. We compute a vector $g\in \R^N$, which is the usual gradient of the functional 
$F(z) = E[t_i, u(z),  m(z),  p(z)] + D [p(z_{i-1}), p(z)]$ 
with the specialty that for $p(z) = p(z_{i-1})$ we select $0$ from $\partial_{z^k} D[p(z_{i-1}), p(z)]$ for $k \in I_p$.
Now, we check if $-g$ actually is a descent direction, that is $F(z-t g)<F(z)$ for sufficiently small $t$. In this case, we perform a descent step in this direction based on a step size controlled line search. Otherwise, we identify the index $k\in I_k$ for which the slope $\lim_{t\to 0} \tfrac{F(z - t g^k) -F(z)}{t}$  is maximal and set $g^k=0$.
Let us remark, that due to our definition of $g$ there is at least one such index with an associated positive slope.
We iterate this until $-g$ is indeed a descent direction and a line search can be performed.
In the line search we prevent an overshooting ensuring that for each $k\in I_p$ the sign of $z^k-z^k_{i-1}$ does not change. 
To this end, if in the line search a selected step size would contradict this property, we replace this step size by the largest step size smaller than the given one 
such that the desired property still holds. This might lead to a new configuration in which the twin boundary does not move, i.e. $p(z) = p(z_{i-1})$.
In our implementation this procedure leads to a robust solution of the minimization problem  \eqref{eq:T}.

\begin{figure}
\includegraphics[width=\linewidth]{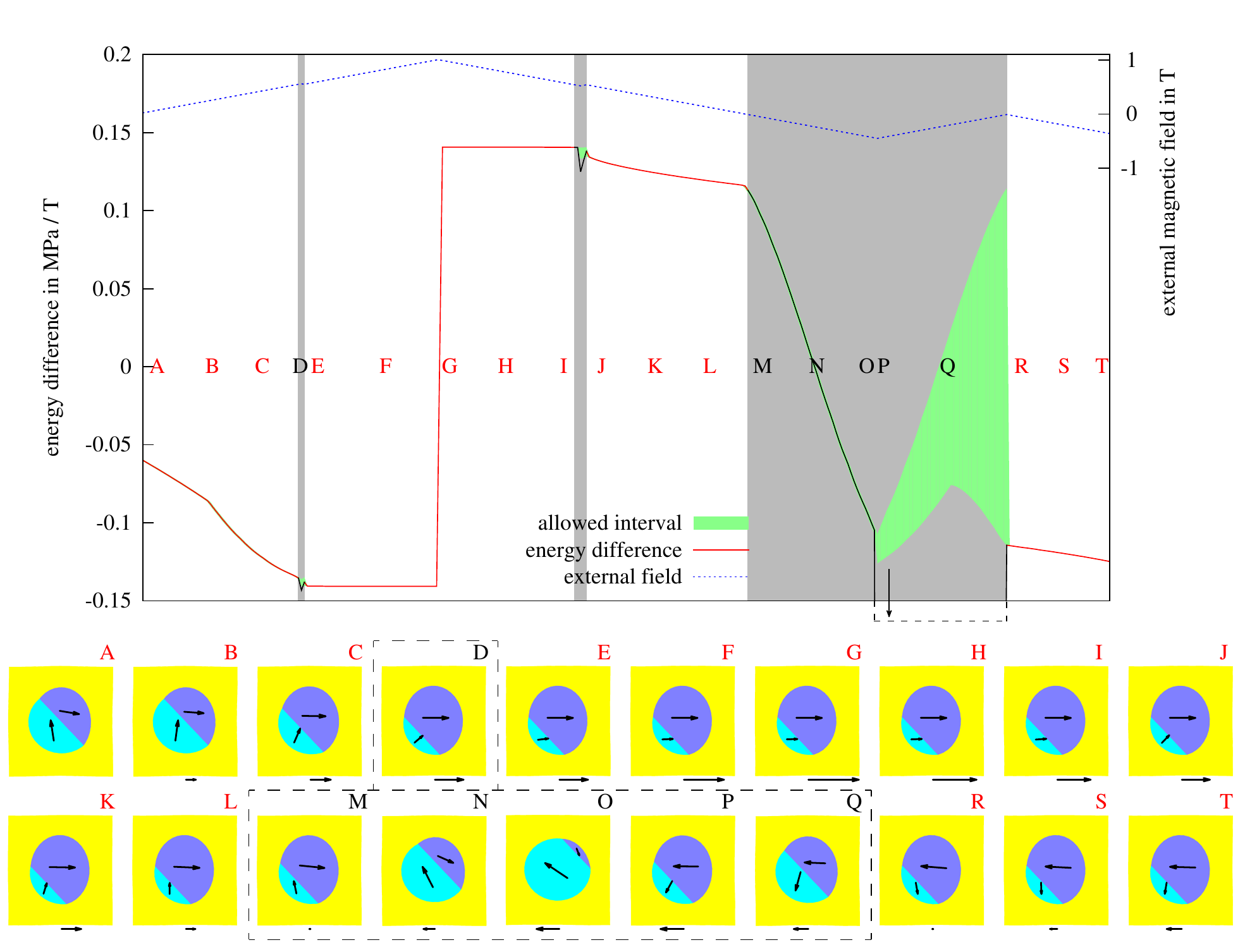}
\caption{The simulation of our rate independent and time discrete hysteresis model which involves a backtracking strategy is displayed. The actual simulation results are depicted below in states (A--C, E--L,R--T) marked in red, whereas the states (D, M--Q) marked in black are later on canceled based on the backtracking. The above diagram shows the energy differences
between consecutive time steps of the simulation divided by the differences in applied field (red/black curve) as well as the applied magnetic field (dotted in blue). The grey regions mark the time intervals effected by the backtracking.}\label{fig:bt}
\end{figure}
\paragraph{Backtracking.}
The incremental minimization problem \eqref{eq:T} demands a global minimization, while the gradient descent algorithm only delivers local minima. In fact, it may happen that the sequence of discrete solutions is continuous along a path of local minima, whose energy deviates substantially from the one associated with the global minimization.
At a later stage the algorithm might jump back into a state of significantly lower energy. This can be detected by the algorithm based on a contradiction to the energy estimates \eqref{eq:EL}  or \eqref{eq:ER}.  In this case one initiates an iteration backward in time. Thereby, one takes into account the new lower energy state as the initial state of the gradient descent scheme in the previous time step aiming for a lower energy state also in that time step. One repeats this iteration step until the energy estimates \eqref{eq:EL} and \eqref{eq:ER} are again fulfilled. Afterwards the time discrete evolution is restarted with the new state at that particular time step. 
This strategy has already been used for related problems in \cite{MielkeRoubicekZeman2010,Benesova2011}. 
 
Let us study this strategy for an example illustrated in Figure \ref {fig:bt}. 
At first, in (A--G) the external magnetic field increases pointing to the right, 
and one observes a growing blue phase (with horizontal easy axis, and magnetization to the right). Next, for decreasing magnetic field (G--M) followed by an increase in the opposite direction (M--O), the algorithm gets stuck in the local minimum with the magnetization of the blue 
phase pointing to the right. Thus, the blue phase is magnetized opposite to the magnetic field and shrinks quickly.  For some critical field (P), the Zeeman 
term is strong enough to pull the magnetization out of its local minimum with a magnetization pointing to the left, now again aligned with the external field. This then leads to an instantaneous growth of the blue phase. 
From the global minimization perspective, switching the 
magnetization to the left would already be favorable at an earlier stage, where the external magnetic is turning from right to left.
The instantaneous growth of the blue phase from (O) to (P) can algorithmically be identified since the energy estimates \eqref{eq:EL} and \eqref{eq:ER} fail. 
In the plot the energy difference between two consecutive time steps is represented by the red/black curve, whereas the green area shows 
for each time step the interval spanned by the different right hand sides of \eqref{eq:EL} and \eqref{eq:ER}.
From (O) to (P) the red/black curve drops down instantaneously leaving the so far thin green interval--as an indication of the above conflict. 
In addition one observed here a spreading of the green area.  
Now we start the backtracking, reported in (P--R). In (R) the energy estimates \eqref{eq:EL} and \eqref{eq:ER} hold again, indicated by the reentry of the red/black curve into the green area. Thus, the actual evolution path is given by removing this wrong forward path (M--O) together with the backtracking (P--Q) (altogether marked in grey) and restarting with the ``better'' local minimum (R).
The final path of the energy difference is plotted in red. The black parts with grey background are local minima that have been improved by backtracking.
Actually, the plot reflects two additional minor occurrences of backtracking at (D), (I--J).

To improve performance, we  try several different starting conditions (e.g. local minima of energy components like the anisotropy and Zeeman energy) 
for the minimization after a couple of timesteps. This permits to find better local minima earlier and therefore avoids long backtracking periods.

\section {Results}
 \label{secresults}
In this section we discuss our numerical results and specific predictions for the behavior of MSM-polymer composites based on the numerical algorithm in two space dimensions presented in Sections \ref{spacediscrete} and \ref{eqc:timediscr}.
Unless otherwise specified, we use
for the polymer an elastic modulus of 
 $E = 1\MPa$, somewhat softer than the polyurethane 
 from \cite{KauffmannScheerbaumLiuetal2012} with $E\sim2\MPa$,
 and a Poisson ratio of $\nu= 0.45$. The periodic cell is a square
 of small side length $\varepsilon$,
 as illustrated in Figure \ref{fig:configuration}, and contains a single
 particle with radius  $0.3\varepsilon$, corresponding to a MSM volume fraction of $28\%$. 
We assume that there are no boundary tractions, in the sense that $g=0$
(see Eq.~(\ref{eqelastmatrix})), and that the 
polymer is stress free if the MSM material is in the austenite phase. This means
that for the polymer the eigenstrain is zero, the two phases of the MSM material are
characterized by the 
eigenstrains in Eq.~\eqref{eqdefeps12}, with suitable rotations; the initial configuration 
is not stress-free.
For the dissipation distance we use
 $D=0.1\MPa$, which corresponds to a switching field of approx $0.25\T$, 
 in agreement with the experimental results of Figure \ref{fig:experiment_simulation}.
We discuss here only the cell problem ignoring the macroscopic energy contributions,
except for Figure \ref{fig:macrostrayfield}. 
The macroscopic strain tensor $A$ entering the affine-periodic boundary conditions
in (\ref{eqbcaffineper}) is assumed to be zero, except for Figure \ref{figure:boundary}.
In most simulations we use a
horizontal magnetic fields up to $1\T$.

We start by discussing Figure \ref{figure:polelastic} in detail, to illustrate
the general features of our results. 
The simulation path is illustrated in the bottom lines of the figure,
with the  initial configuration (A) on the left. 
The MSM particle is divided into two phases, with different magnetizations (shown by the arrows)
and eigendeformations (in the reference configuration the particle is a disc).
In this initial configuration, without external field, the two variants are
equivalent by symmetry. The anisotropy energy would be minimized by a
horizontal and vertical orientation of the magnetization in the two variants.
This configuration would not generate any magnetic charge at the twin boundary, since
the normal component would be continuous. However, there would be significant magnetic charges
at the boundaries between the MSM and the polymer. These, and correspondingly the demagnetization
energy $E_{\R^2}^\demag$, are reduced by 
rotating the magnetization towards the interface.

Starting from the initial configuration (A) we apply an horizontal magnetic field, which increases up to $1\T$.
The presence of an external fields favors one of the variants, which has
an almost horizontal magnetization. Already for small fields  the
phase boundary moves (B). At the same time the magnetization inside each
phase rotates (B--C), to better accommodate the external field. 
Correspondingly the fraction of the first variant (which is 
closely related to the eigenstrain) increases, see red curve in the first plot, and the horizontal magnetic moment also increases, see 
red curve in the right plot. For fields
close to $1\T$ the minority phase has almost completely disappeared (C).

Reducing the field the particle only transforms back with some delay (D). 
After the field is removed (E), the two phases are again equivalent from the viewpoint of anisotropy. The purely magnetic material 
would have no reason to transform back. However, in the composite
the polymer elasticity favors the initial state, in which the average strain is zero.
This is contrasted by dissipation, which resists any movement of the interface. 
Thus the interface is only partially pushed back at zero field
(E); we discuss below a situation in which the field changes orientation. Increasing now the field in the opposite direction (E--G)
the magnetization in the majority phase flips to the other sign; the same
phase is favored for all horizontal fields, independently on the orientation.
Therefore the volume of the ``horizontal'' phase increases again to almost all of the particle (G). Decreasing the field again to zero shows the same behavior, 
and iteration gives a repetition of the
hysteresis loop (C--G), without ever going back to the initial state (A).

\bigskip\noindent\begin{minipage}{\linewidth}\centering
\includegraphics[width=0.8\linewidth]{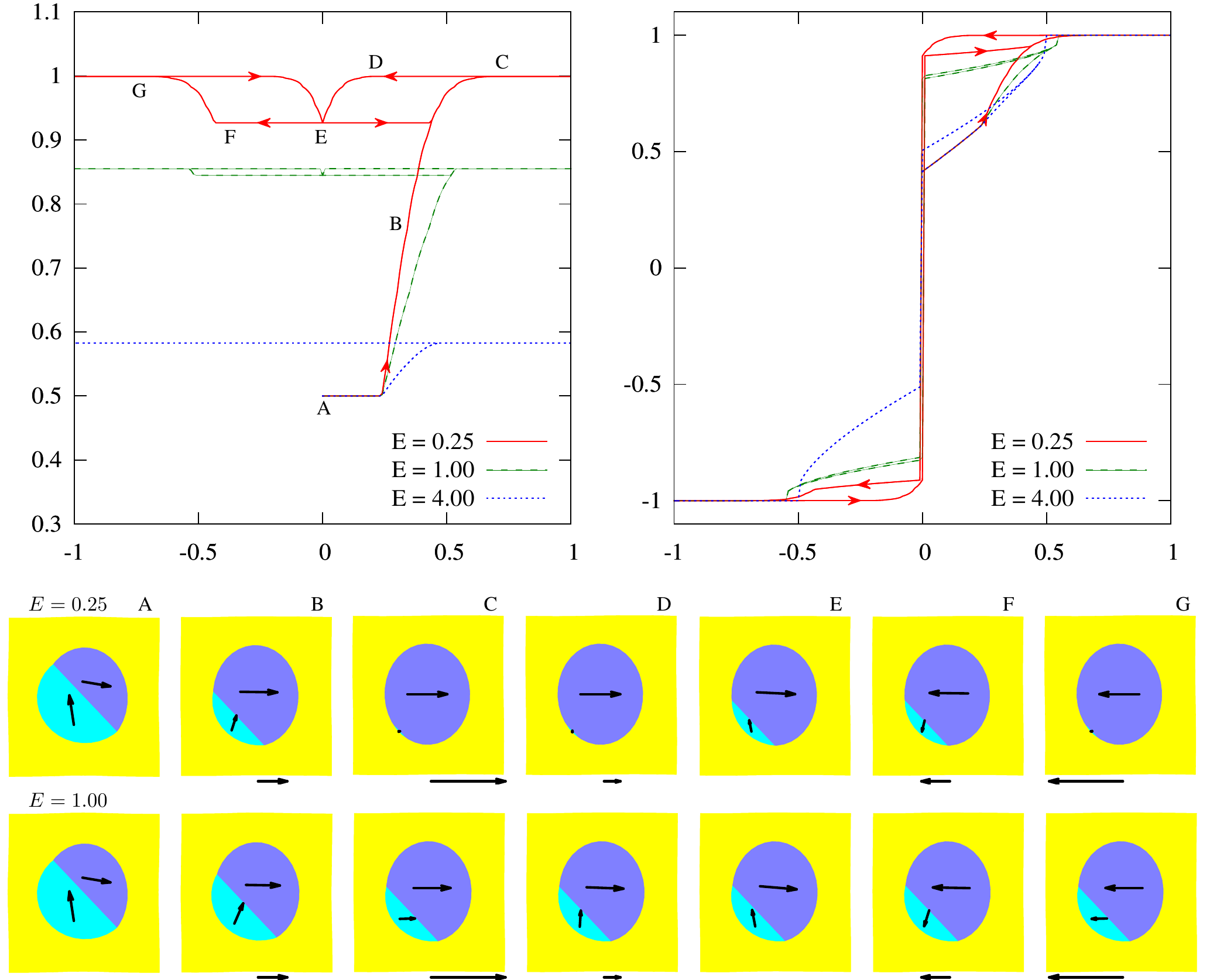}
\captionof{figure}{Simulation of the hysteresis loop (states C--G), including an initial phase (states A,B), in a MSM-polymer composite 
with one disc shaped particle per cell of the periodic lattice. 
Parameters as given at the beginning
of Section \ref{secresults}. For the polymer
elasticity modulus we compare three different values, the standard one $E=1\MPa$,
a larger value $E=4\MPa$ and a smaller value $E=0.25\MPa$.
The plot on the left side depicts the volume fraction of one variant over the strength $H(t)$ of the external magnetic field in $\T$. The plot on the right side shows the horizontal component of the average magnetization (relative to saturation) again over $H(t)$.}
\label{figure:polelastic}
\end{minipage}\bigskip

\bigskip\noindent\begin{minipage}{\linewidth}\centering
\includegraphics[width=0.7\linewidth]{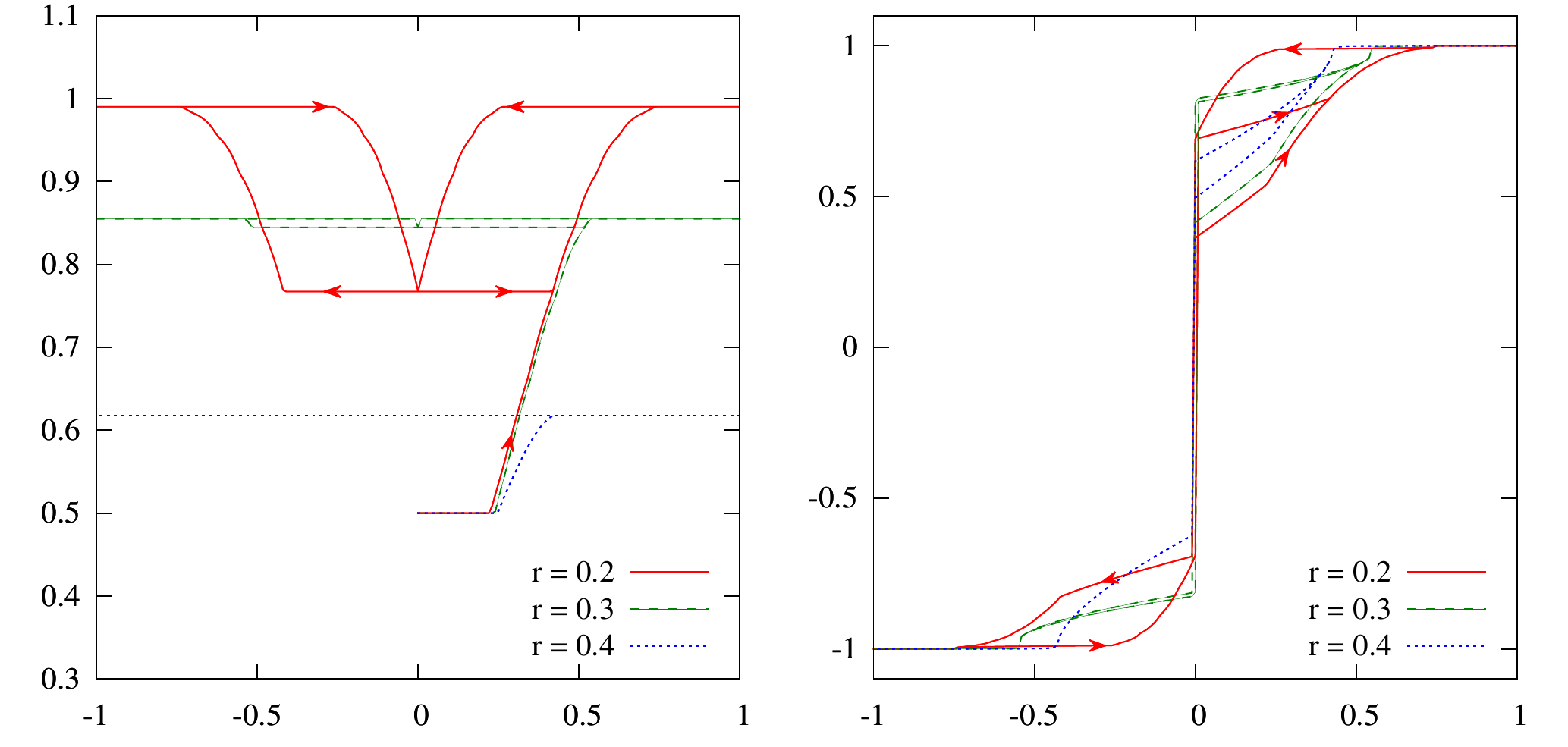}
\captionof{figure}{Simulations with different particle radius.
Same parameters and geometry as in Figure \ref{figure:polelastic}. The left plot again shows the volume fraction of one variant, the right plot the horizontal part of the magnetization.}
\label{figure:radius}
\end{minipage}\bigskip

In the following we illustrate the effect of the various parameters, which permits to better understand
the mechanisms behind the observed cycle, and to make the role of the different terms quantitative.  The general picture in most of our simulations is similar to the one
we just described, therefore we focus on the resulting diagrams and highlight the differences.

In Figure \ref{figure:polelastic} we also show curves obtained for different 
values of the polymer elasticity. The amplitude of the hysteresis loop is significantly
reduced for stiffer polymer. This corresponds to the fact that the magnetic anisotropy
is not sufficient to push the interface far away from the initial position. The 
magnetization however becomes fully aligned with the field in both variants, as the 
right plot shows.

Figure \ref{figure:radius} shows the influence of the particle radius, which corresponds to the MSM volume fraction, on the hysteresis loop. The amplitude of the deformation is largest for the 
smallest volume fraction. Indeed, since no macroscopic deformation is possible, the motion of the interface would be
almost completely inhibited if the particle would fill the entire simulation cell. At the same time, the largest particles
generate a larger work output, as was made quantitative in \cite{ContiLenzRumpf2007}.

\bigskip\noindent\begin{minipage}{\linewidth}\centering
\includegraphics[width=0.7\linewidth]{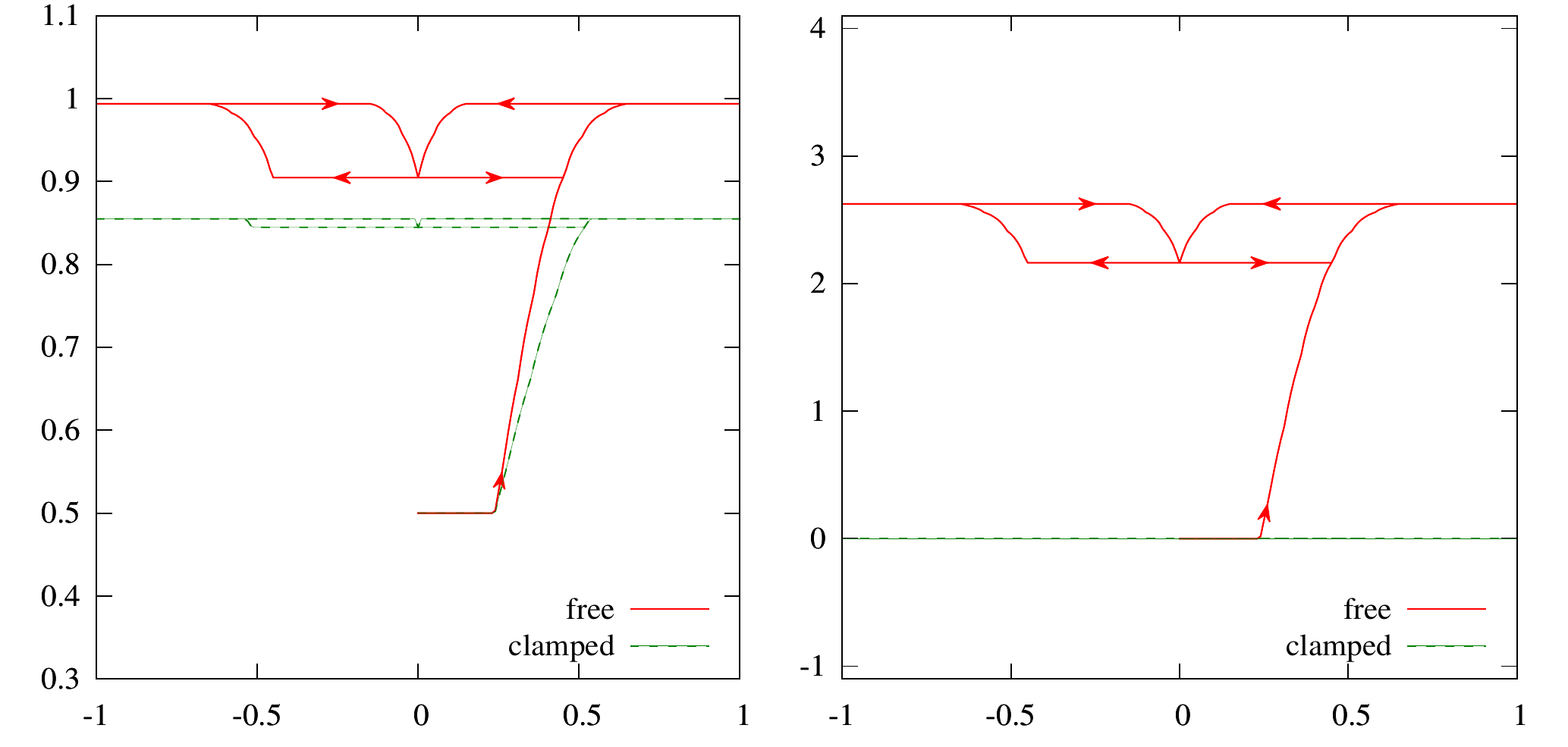}
\captionof{figure}{Simulations with different macroscopic elastic boundary conditions.  The case
 $A=0$ is compared to the case, where the energy is also minimized with respect to the affine macroscopic strain tensor $A$ (same parameters and geometry as in Figure \ref{figure:polelastic}).
On the left we show the volume fraction, on the right the macroscopic deformation.}
\label{figure:boundary}
\end{minipage}\bigskip

In Figure \ref{figure:boundary}  we investigate the role of boundary conditions. Following the homogenization paradigm we solve the cell problem 
 with affine-periodic boundary data, as given in (\ref{eqbcaffineper}). The affine matrix  $A\in \R^{2\times 2}$ corresponds 
to the local macroscopic strain tensor. In most of our simulations  $A$ is zero. 
Here we compare with a simulation in which we minimize also over the matrix $A$, leading to a
spontaneous strain of $2.6\%$.

\begin{wrapfigure}{r}{0.35\linewidth}\ \\[-5mm]
\includegraphics[width=\linewidth]{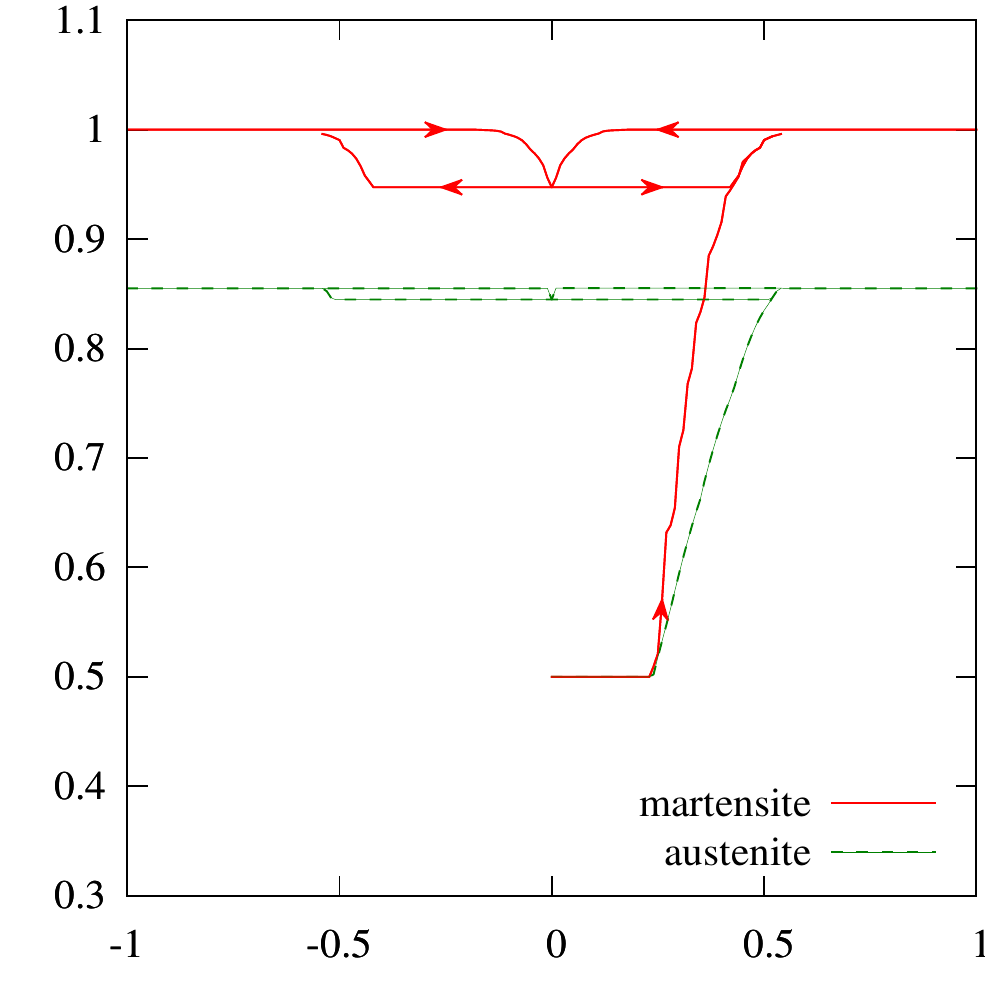}
\captionof{figure}{Simulations with different elastic reference configurations.
Same parameters and geometry as in Figure \ref{figure:polelastic}, plot of volume fraction.}
\label{figure:reference}
\end{wrapfigure}
As discussed in the beginning of this section, we usually assume the polymer to be stress free in the austenite phase. The simulation is always performed for the martensite phase, thus even in the initial configuration the polymer is not stress free. In Figure  \ref{figure:reference} we compare this to the case that the polymer is stress free in the martensite phase, i.e. there is no pre-stress on the polymer. Without pre-stress the polymer is easier to deform, and thus a larger fraction of the particle is actually transformed. 

\begin{wrapfigure}{l}{0.35\linewidth}\ \\[-5mm]
\includegraphics[width=\linewidth]{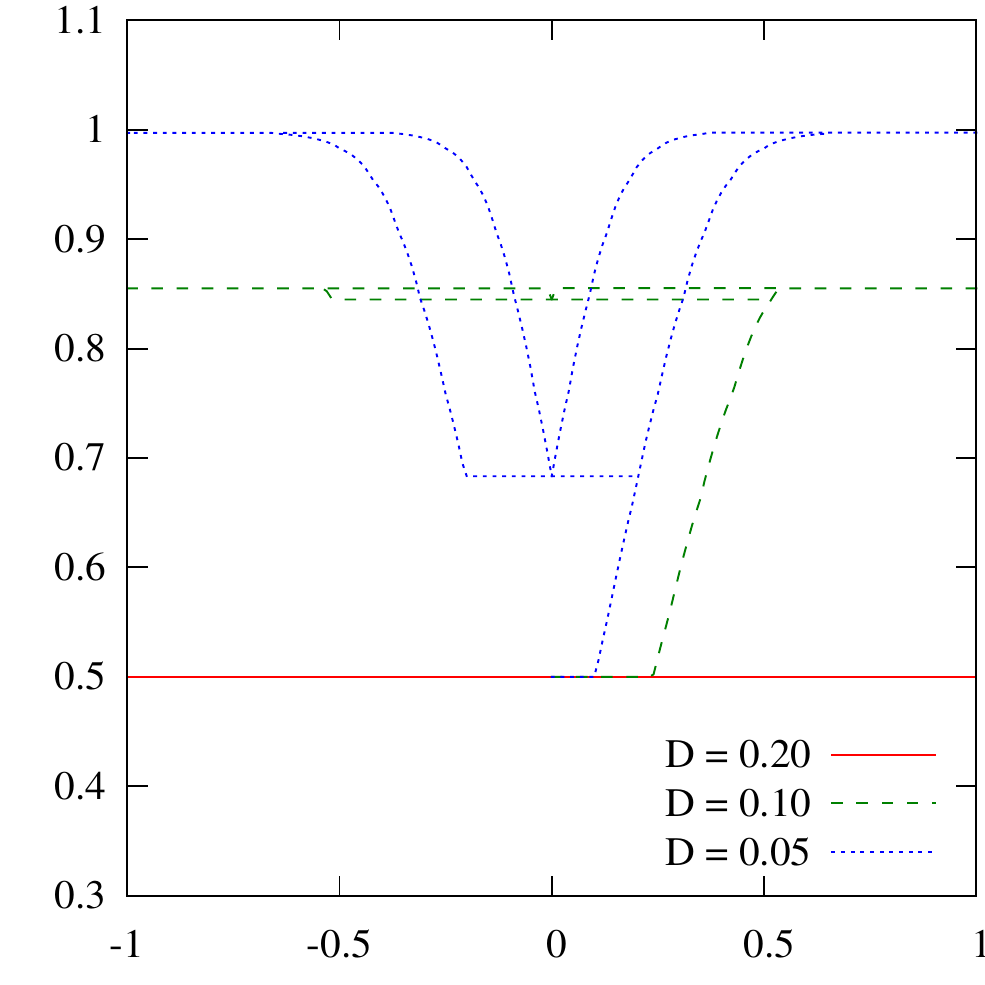}
\captionof{figure}{For different dissipation coefficients the volume fraction is plotted (parameters and geometry as in Figure \ref{figure:polelastic}, $E=1\MPa$).}
\label{figure:dissipation}
\end{wrapfigure}

Figure \ref{figure:dissipation} shows that the strength of the dissipation also has a large
effect on the hysteresis. For $D=0.20\MPa$ the transformation is blocked, the domain wall 
does not move from the initial position. This corresponds to the fact that the magnetic anisotropy  with $K_u=0.13\MPa$ is smaller than the dissipation. 

Figure \ref{figure:horientation} shows the effect of the orientation of the applied field, which allows 
to characterize the alignment between the lattice orientation and the macroscopic applied field. Here, the particle configuration is fixed with horizontal and vertical easy axes, and the external magnetic field
increases in a fixed direction, that is rotated by some angle. If the rotation is $5°$, a larger field is required to start the transformation, but the general behavior is otherwise very similar. For a larger angle of 
$10°$, the transformation starts even later and is much slower. Indeed, in this case, the polymer is not strong enough to initiate a significant backwards transformation, so there is no repeatable hysteresis when the field changes to the opposite direction. For an even larger rotation of $20°$, there is no observable movement at all.

In Figure \ref{fig:macrostrayfield} we consider the stray field including the macroscopic part, for the case of a circular composite workpiece.
As discussed in \cite{ContiLenzRumpf2007}, the macroscopic part of the stray field can be effectively computed in this case. It can be interpreted as a modification of the effective external field in the cell problem,
which reduces the strength of the magnetic field that actually acts on the particle. Because of this a larger external field ($1.5\T$ instead of $0.6\T$) is needed to reach magnetic saturation, and the transformation happens
much slower. The backward transformation, however, in fact starts earlier and goes further compared to the simulation without macroscopic stray field. 
This is mainly due to the fact that we do not consider solutions with larger periodicity or magnetic domain walls within the different phases. 
Because of this, a pattern of alternating magnetizations (cf. \cite{ContiLenzRumpf2007}) that significantly reduces the macroscopic stray field energy can only be realized by transforming back 
towards a 50:50 distribution of the phases and their associated magnetizations.

Finally, in Figure \ref{figure:horizvertic} we illustrate the results of simulations in which the field
is biaxial. We start from field zero, then increase it up to $1\T$ horizontally, afterwards bring it back to zero,
and finally  increase it up to $1\T$ in the vertical direction. Whereas the field favors at first the
``horizontal'' martensitic variant, in the second round the vertical field favors the other variant. Therefore
the back transformation of the interface is not only due to the polymer, but also due to the Zeeman term in the energy. 
This leads to a substantially larger hysteresis cycle.

\bigskip\noindent\begin{minipage}{\linewidth}\centering
\includegraphics[width=0.7\linewidth]{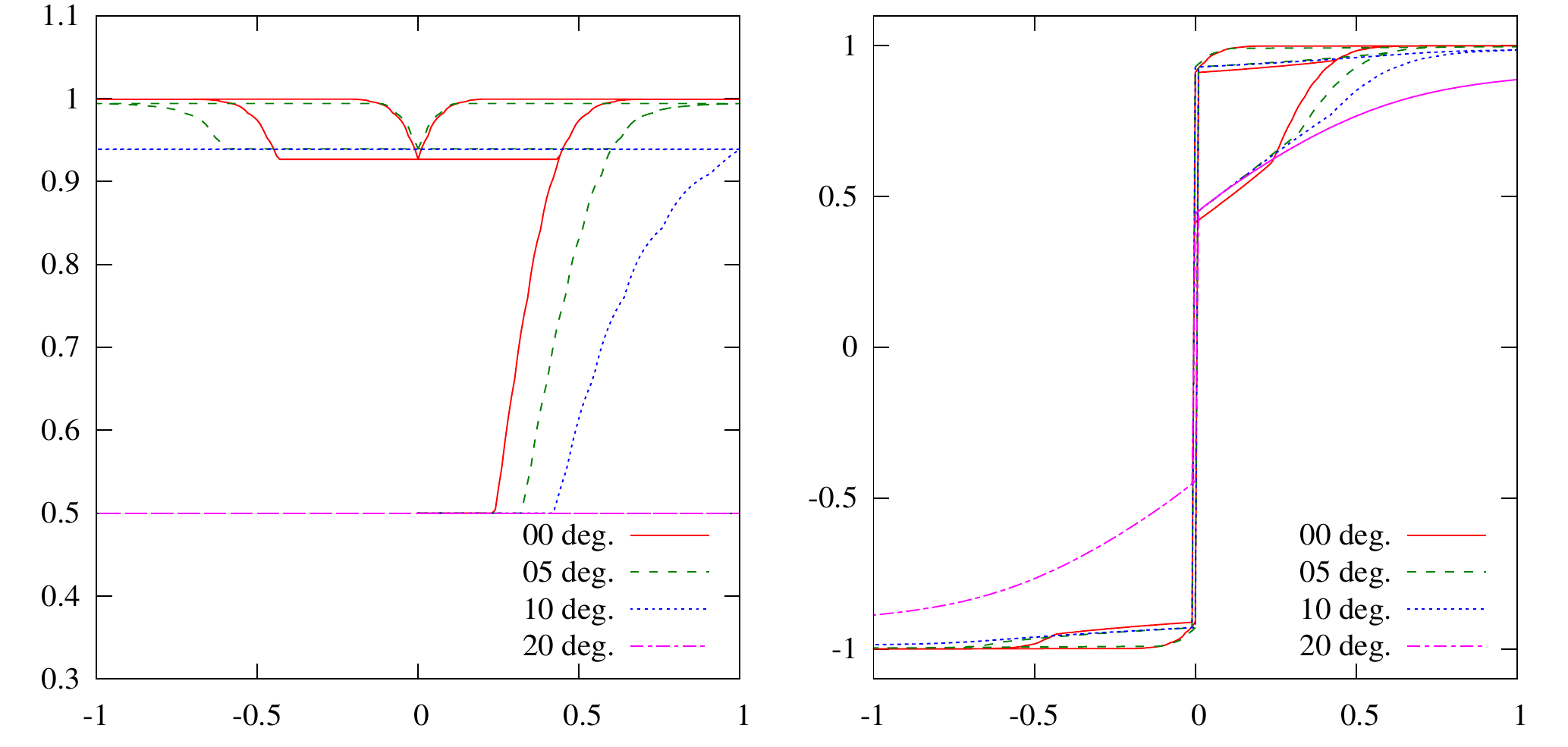}
\captionof{figure}{Simulations with different degrees of misalignment between the external magnetic field and the magnetic easy axis of one phase, those easy axis is horizontal. We plot the volume fraction of one variant (left) and the horizontal component of the magnetization (right) for different rotations of the external field.
Same parameters and geometry as in Figure \ref{figure:polelastic}.}
\label{figure:horientation}
\end{minipage}\bigskip

\bigskip\noindent\begin{minipage}{\linewidth}\centering
\includegraphics[width=0.7\linewidth]{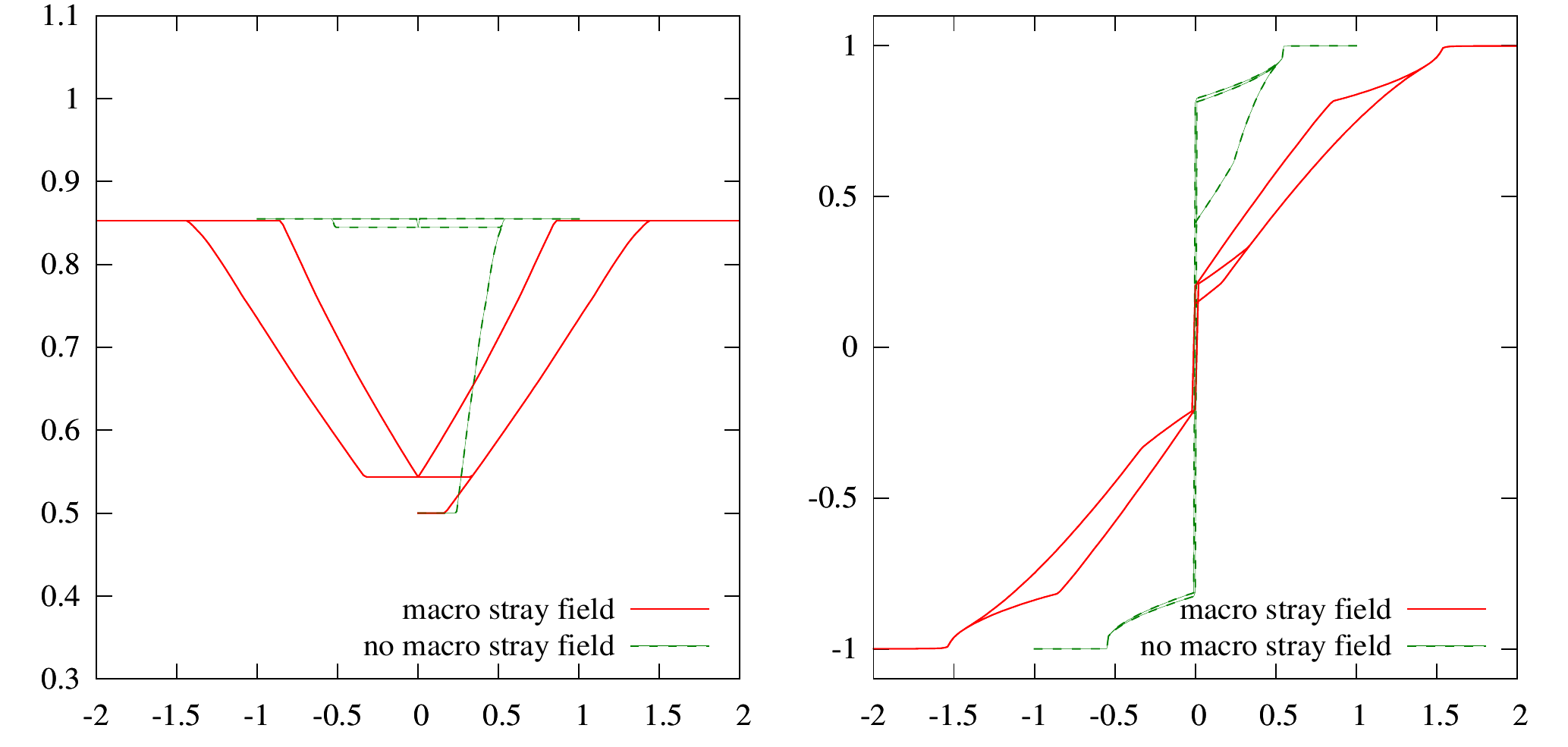}
\captionof{figure}{Influence of the macroscopic stray field is shown plotting the volume fraction of one variant on the left and the horizontal component of the magnetization on the right.
Same parameters and geometry as in Figure \ref{figure:polelastic}.}
\label{fig:macrostrayfield}
\end{minipage}\bigskip

\bigskip\noindent\begin{minipage}{\linewidth}\centering
\includegraphics[width=0.7\linewidth]{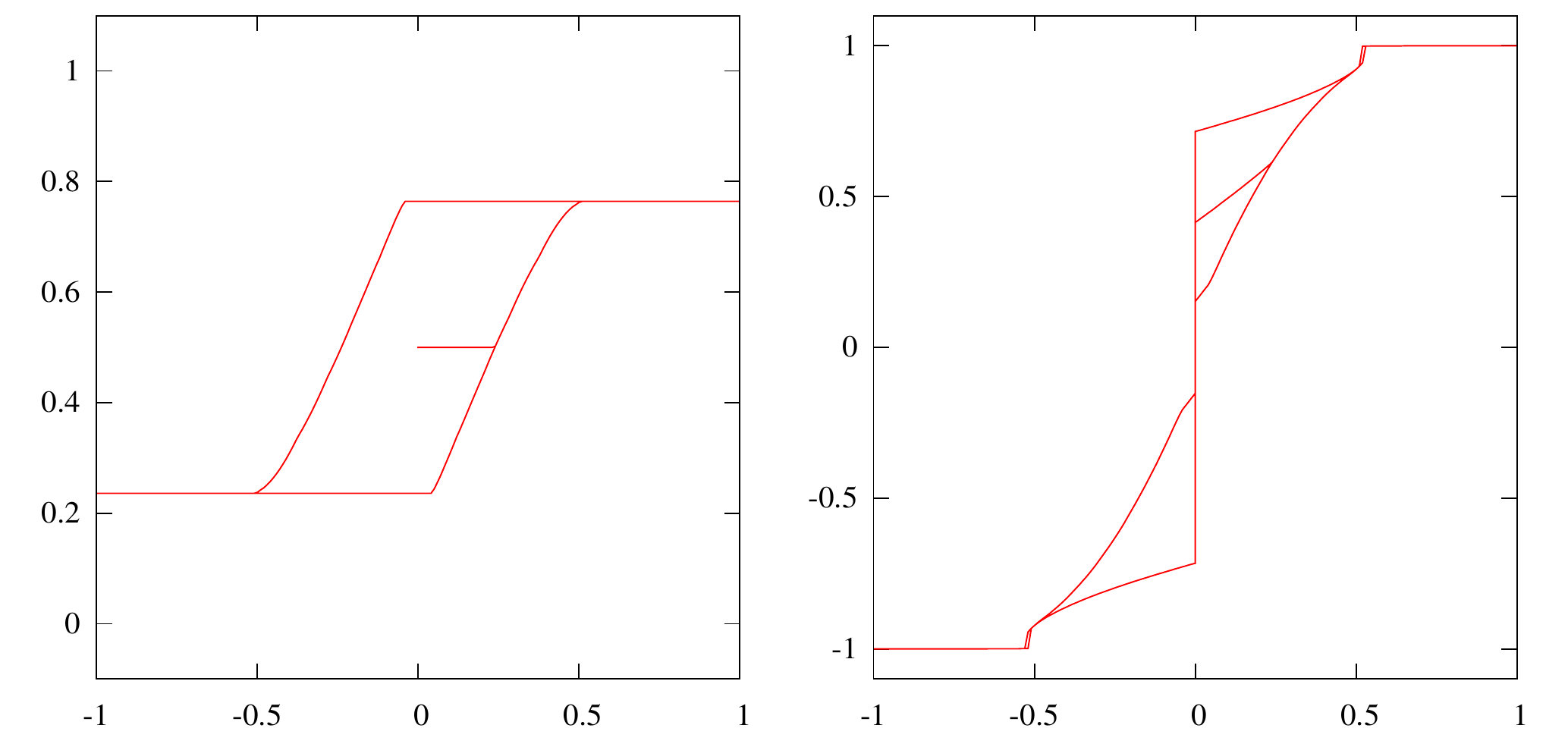}
\captionof{figure}{Simulations with external field first along the $x$ axis and then along the $y$ axis 
plotting the volume fraction of one variant on the left and the horizontal component of the magnetization on the right. The negative
values of the external field and of the magnetization in both plots represent fields and magnetizations along the $y$ axis. 
Same parameters and geometry as in Figure \ref{figure:polelastic}.}
\label{figure:horizvertic}
\end{minipage}\bigskip

\section{Conclusions}
We have presented and studied a rate-independent model for MSM-composite materials.
Our results show that phase transformation and hysteresis are largely influenced by the material and geometric parameters. In particular, the transformation is inhibited by large polymer elasticity 
coefficients and by large dissipation. Furthermore, the deformation is enhanced by using an experimental protocol in which external fields in two orthogonal directions are used. In contrast, 
small misalignment between the field and the particle lattice orientation seems to play a minor role. Although our simulations are restricted to two spatial dimensions, we expect the general results and the trends we identified to be valid also in three dimensions. 
We hope that our findings may be helpful in the experimental search for 
composite materials with large spontaneous strains.

\section*{Acknowledgments}
This work was partially supported by the Deutsche Forschungsgemeinschaft
through  the Schwerpunktprogramm 1239 {\em \"Anderung von
  Mikrostruktur und Form fester Werkstoffe durch \"au{\ss}ere Magnetfelder}.

\bibliography{msm,bibtex/all,bibtex/own}
\bibliographystyle{alpha-noname-nonumber}

\end{document}